\documentclass{amsart}
\usepackage{graphicx}
\theoremstyle{definition}

\theoremstyle{remark}

\numberwithin{equation}{section}

\begin{document}

\title{ON DIRECTIONAL ENTROPY OF A $\mathbb{Z}^{2}$-ACTION}%
\author{Hasan Akin}%
\address{}%
\email{}%

\thanks{}%
\renewcommand{\thefootnote}{}


\begin{abstract}

Consider the cellular automata (CA) of $\mathbb{Z}^{2}$-action
$\Phi$ on the space of all doubly infinite sequences with values
in a finite set $\mathbb{Z}_{r}$, $r \geq 2$ determined by
cellular automata $T_{F[-k, k]}$ with an additive automaton rule \\
$F(x_{n-k},...,x_{n+k})=\sum\limits_{i=-k}^{k}a_{i}x_{n+i}(mod
r)$. It is investigated the concept of the measure theoretic
directional entropy per unit of length in the direction
$\omega_{0}$. It is shown that
$h_{\mu}(T_{F[-k,k]}^{u})=uh_{\mu}(T_{F[-k,k]})$,
$h_{\mu}(\Phi^{u})=uh_{\mu}(\Phi)$ and
$h_{\vec{v}}(\Phi^{u})=uh_{\vec{v}}(\Phi)$ for $\vec{v} \in
\mathbb{Z}^{2}$ where $h$ is the measure-theoretic entropy.

\end{abstract}
\maketitle

\section{Introduction}

\par In the present paper we study directional entropy of
$\mathbb{Z}^{2}$-action generated by an additive cellular automata
(CA). CA initialed by Ulam and von Neumann has been investigated
by Hedlund [4]. He systematically studied purely mathematical
point of view. In Hedlund's work are given current problems of
symbolic dynamics. In [7], Shereshevsky has investigated ergodic
properties of CA, and also defined the n-th iteration of a
permutative cellular automata.
\par The concept of the directional entropy of a $\mathbb{Z}^{2}$-action
has first been introduced by Milnor [5]. Milnor defined the
concept of the directional entropy function for
$\mathbb{Z}^{2}$-action generated by a full shift and a block map.
This concept was also studied in [2], [6] and [8].
\par In [2], Courbage and Kaminski have calculated the directional
entropy for any cellular automata (CA) of $\mathbb{Z}^{2}$-action
$\Phi$ on the space of all doubly infinite sequences with values
in a finite set $A$, determined by an automaton rule $F[l,r]$,
$l,r\in \mathbb{Z}$, $l\le r$, and any $\Phi$-invariant Borel
probability measure.  In [6], Park expressed the directional
entropy in an integral form.
\par In [1], the author calculated the measure entropy of
additive one-dimensional cellular automata with respect to uniform
Bernoulli measure. In [3], Coven and Paul investigated some
properties of the endomorphisms of irreducible subshifts of finite
type and n-block maps.
\par The shift $\sigma$ and $T_{F[-k,k]}$
are commuted and if $T_{F[-k,k]}$ is non-invertible, they generate
a $\mathbb{Z}\times\mathbb{N}$ action $\Phi^{(p,q)}=\sigma^{p}
T^{q}_{F[-k,k]}$, which can be extended to $\mathbb{Z}^{2}$-action
on $\Omega$. Notice that
$\sigma^{i}T_{F[-k,k]}=T_{F[-k,k]}\sigma^{i}=T_{F[-k+i,k+i]}$ for
all $i\in \mathbb{Z}$. We suppose that $\mu$ is a probability
ergodic measure which is invariant under the action $\Phi$ of
$\mathbb{Z}\times\mathbb{N}$. Let $\vec{v}$ be an arbitrary vector
of $\mathbb{Z}^{2}$. Denote by $h_{\vec{v}}(\Phi)$ the directional
entropy of $\Phi$ [2]. The measure-theoretic entropy of $\Phi^{(p,
q)}$ with respect to $\mu$ is denoted by $h_{p,q} =
h(\sigma^{p}T_{F[-k,k]}^{q} )$ where $F^{n}$ denotes the n-th
iteration of a function (or map) $F$ (cf. [7]). It is easy to show
that $h_{p,0}<\infty$ for all $-\infty <p <\infty$.
\par The question posed by Milnor [5]
is :  Does the limit
$$\lim\limits_{i \to
\infty}\frac{1}{\sqrt{m_{i}^{2}+n_{i}^{2}}}h_{m_{i},n_{i}}
$$
exist for the sequence $\{(m_{i},n_{i}), \, i=\overline{1,
\infty}\}\subset \mathbb{Z}^{2}$, $m_{i}\rightarrow \infty $,
$n_{i}\rightarrow \infty$, $\frac{m_{i}}{n_{i}}\rightarrow \omega
_{0}$ as
 $i\rightarrow \infty $?  \par An affirmative answer to this question
was given by Park [6] and Sinai [8] for an irrational number
$\omega _{0}$. Sinai [8] and Park [6] also showed that the
function $h_{p,q}$ is a homogeneous function of the first degree,
i.e. $h_{up,uq} =|u|h_{p,q}$.

\par In this paper under additional assumptions we show that
$h_{\mu}(T_{F[-k,k]}^{u})=uh_{\mu}(T_{F[-k,k]})$,
$h_{\mu}(\Phi^{u})=uh_{\mu}(\Phi)$ and
$h_{\vec{v}}(\Phi^{u})=uh_{\vec{v}}(\Phi)$.

\section{Preliminaries}
Let $\mathbb{Z}_{r}=\{0, 1,..., r-1\}$ be the set of integers
modulo $r$ and denotes a state set of each cell and
$\Omega=\prod\limits_{i=-
\infty}^{\infty}\mathbb{Z}_{r}=\mathbb{Z}_{r}^{Z}$ be the space of
all doubly infinite sequences\\
$x=\{x_{i}\}_{i=-\infty }^{\infty }$, $x_{i} \in \mathbb{Z}_{r}$.
$\Omega$ is compact in topology of direct product and a measurable
space. We denote by $\sigma$ the shift transformation on $\Omega$,
i.e. $(\sigma x)_{i}=x_{i+1}$ for $i \in \mathbb{Z}$.
\par It is obvious that  $\sigma$ is homeomorphism.
Let $\textbf{M}$ be the product $\sigma$-algebra of $\Omega$ and
$\mu$ be a probability invariant measure. The quadruplet
($\Omega$, $\textbf{M}$, $\mu$, $\sigma$) is called symbolic
dynamic system.

\par Let $m$ be a fixed positive integer. We denote by
$\mathbb{Z}^{m}_{r}$ the $m$-fold direct product $\mathbb{Z}_{r}
\times\ldots\times\mathbb{Z}_{r}$.
\par An automaton rule $F$ is said to be right permutative (cf. [7])
if for any\\
$(\bar{x}_{1},\ldots,\bar{x}_{m-1}) \in \mathbb{Z}^{m-1}_{r}$ the
mapping $x_{m} \rightarrow
F(\bar{x}_{1},\ldots,\bar{x}_{m-1},x_{m})$ is a permutation of
$\mathbb{Z}_{r}$. Similarly we define a left permutative mapping.
\par We say that $F$ is bipermutative if it is right and left permutative.
\par Any mapping $F: \mathbb{Z}^{2k+1}_{r} \to
\mathbb{Z}_{r}$ is called an automaton rule. Take any nonnegative
integer $k$ and consider a linear map $F:\mathbb{Z}^{2k+1}_{r} \to
\mathbb{Z}_{r}$ defined by formula
$$
F(z_{-k},...,z_{k})=\sum\limits_{i=-k}^{k} a_{i}z_{i}({mod} r)
\eqno(1)
$$
where $a_{i}\in \mathbb{Z}_{r}$, $i=\overline{(-k, k)}$. An
automaton rule $F$ in the form (1) is called an additive automaton
rule.
\par The homeomorphism $T_{F[-k,k]}: \Omega \to \Omega$
defined as
$$
(T_{F[-k,k]}x)_{n}=F(x_{n-k},...,x_{n+k}) = \sum\limits_{i=-k}^{k}
a_{i}x_{n+i} (mod r), n\in \mathbb{Z}
$$
is said to be the additive one-dimensional cellular automata (CA)
defined by $F[-k, k]$.
\par It it clear that the additive CA-map $T_{F[-k,k]}$ is surjective and
non-invertible. Moreover, $T_{F[-k,k]}$ preserves the uniform
Bernoulli measure $\mu$ [7].
\par In [7], Shereshevsky has define inductively
the u-th iteration
$F^{u}:\mathbb{Z}_{r}^{2ku+1}\rightarrow\mathbb{Z}_{r}$ of the
rule $F$ as follows:

\begin{flushleft}
$F^{u}(x_{-2ku},...,x_{-2ku+2k},...,x_{-2ku+4k},...,x_{2ku-2k},...,x_{2ku})=
F^{u-1}(F(x_{-2ku},...,x_{-2ku+2k}),F(x_{-2ku+1},...,x_{-2ku+2k+1}),...,
F(x_{2ku-2k},...,x_{2ku})).$
\end{flushleft}

\par \textbf{Lemma 2.1.} ([7], Lemma 1.6) The u-th iteration
$T_{F[-k,k]}^{u}$ of CA-map $T_{F[-k,k]}$ generated by the
rule $F$ coincides with the CA-map $T_{F^{u}[-ku,ku]}.$ \\
\par It can be easily checked that the shift
$\sigma$ and a cellular automaton map $T_{F_{[-k,k]}}$ are
commutted i.e. $\sigma\circ T_{F[-k,k]}=T_{F[-k,k]}\circ \sigma$.
The $\mathbb{Z}^{2}$-action $\Phi$ generated by $\sigma$ and
$T_{F[-k,k]}$, i.e. $\Phi^{(p,q)}=\sigma^{p}T_{F[-k,k]}^{q}$ is
said to be a CA-action, if $T_{F[-k,k]}$ is invertible.
\par Let ($\Omega$, $\textbf{M}$, $\mu$, $\sigma$) be a
symbolic dynamic system. Let $\prec$ denotes the lexicograpical
ordering of $\mathbb{Z}^{2}$. Denote by $O$ the zero of
$\mathbb{Z}^{2}$. A sub $\sigma$-algebra $\textbf{A}$ is said to
be invariant if $\Phi^{(p,q)}(\textbf{A})\subset \textbf{A}$ for
every $(p,q) \prec O$. It is clear that $\textbf{A}$ is invariant
iff $\sigma^{-1}(\textbf{A})\subset \textbf{A}$ and
$T_{F[-k,k]}^{-1} (\textbf{A}) \subset \textbf{A}$.
\par Let $\xi$ be a zero-time partition of $\Omega$;

$$
\xi=\{C_{0}(0), C_{0}(1), \ldots, C_{0}(r-1) \}
$$
where $C_{0}(i)=\{x\in \Omega; \, x_{0}=i \}$, $i \in
\mathbb{Z}_{r}$, is a cylinder set.
\par We note that if cellular
automata $T_{F[-k,k]}$ is permutative then the partition\\
$\xi=\{C_{0}(i), \, i\in \mathbb{Z}_{r} \}$ is a generating
partition for CA-map $T_{F[-k,k]}$.
\par Now we introduce some necessary notations. Let $a \in R^{1}$,
$\omega \in R^{+}$, and \\
$I=I(a, \omega)$ be a closed interval on the plane with endpoints
$(a, 0)$ and $(a+{\omega}^{-1}, 1)$, and $\Gamma(a, \omega)$ be a
half-line $y=\omega(x-a)$, $y \le 1$. Suppose that a probability
measure $\mu$ on $\mathbf{M}$ is invariant with respect to the
shift $\sigma$ and cellular automata $T_{F[-k,k]}$.
\par Define the following conditional properties:
$$
 H_{r}(I)=H(
\bigvee\limits_{a+\omega ^{-1}\leq p} \Phi ^{(p,1)}\ \xi\ \mid\
\bigvee\limits_{q=0}^{\infty} \,\bigvee\limits_{a+\omega^{-1}q \le
p} \Phi ^{(p,-q)}\xi)
$$
$$
H_{l}(I)=H( \bigvee\limits_{p \le a+\omega ^{-1}}\Phi ^{(p,1)}\
\xi\ \mid\ \, \bigvee\limits_{q=0}^{\infty} \, \bigvee\limits_{p
\le a+\omega^{-1}q} \Phi ^{(p,-q)}\xi )
$$
where $H_{r}(I)$ and $H_{l}(I)$ are called the right and left
entropies, respectively. In [8], it was shown that these entropies
are finite.
\par Let $(p, q)$ be a point of $\mathbb{Z}^{2}$. Denote $h_{p,
q}=h({\Phi}^{(p, q)})=h({\sigma}^{p}T^{q}_{F[-k,k]})$. The value
of $h_{p, q}$ is equal to the limit
$$
h_{p, q}=\lim\limits_{s \to \infty} H_{\mu} \left(
\bigvee\limits_{n=1}^{q} \, \bigvee\limits_{|m-(a+\omega
^{-1}n)|\leq s} \Phi ^{(m,n)}\ \xi \ \mid\
\bigvee\limits_{n=0}^{\infty} \, \bigvee\limits_{|m-(a+\omega
^{-1}n)|\leq s} \Phi ^{(m,-n)}\xi \right)
$$
with $\omega=\frac{p}{q}$.
\par Let ${\omega}_{0}$ be an
irrational number, $\{(m_{i}, n_{i}), \, i=\overline{0, \infty}
\}$be a sequence of points of the lattice $\mathbb{Z}\times
\mathbb{N}$ such that $m_{i} \to + \infty$ or $m_{i} \to -
\infty$, $n_{i}\to + \infty$ and $\lim\limits_{i \to \infty}
\frac{m_i}{n_i}={\omega} _{0}$.
\par Sinai has proved in [8] that there exists a finite limit
$$
\lim\limits_{i \to
\infty}\frac{1}{\sqrt{m_{i}^{2}+n_{i}^{2}}}h_{m_{i},n_{i}}=C
\eqno(2)
$$
and it doesn't depend on the choise of the sequence $\{(m_{i},
n_{i}) \}.$
\par \textbf{Definition 2.2.} The value $C$ of the limit (2) is
called an entropy per unit of length in the direction
${\omega}_{0}$.
\par It is well known that the automaton map is not one-to-one, in
general, so we should consider the natural extension of the
automaton map (determined by an automaton rule), we need to use
the natural extension the semi-group action to a group action.
\par Let $(\hat{\Omega},\hat{\textbf{M}},\hat{\mu},\hat{T})$
be a natural extension of the dynamical system
$(\Omega,\textbf{M},\mu,T)$ (cf. [2])

Let us recall that $\hat{T}$ is defined as follows:
$$
\hat{T}\hat{x} = (Tx^{(0)},x^{(0)},\ldots), \hat{x} =
(x^{(0)},x^{(1)},\ldots)
$$
where $Tx^{(i)} = x^{(i-1)}, i\geq 1.$ We put
$$ \hat{\tau}\hat{x} = (\tau x^{(0)},\tau x^{(1)},\ldots).
$$
Obviously, $\hat{\tau}\hat{T} = \hat{T}\hat{\tau}$. The
$\mathbb{Z}^{2}$ - action $\Phi$ generated by $\hat{\tau}$ and
$\hat{T}$:
$$
\Phi^{(p,q)} = \hat{\tau}^{p}\hat{T}^{q}
$$
is said to be a CA-action. For a positive integer $m$ and
$E\in\textbf{M}$ we put
$$
E^{(m)} = \{\hat{x}\in \hat{\Omega}; x^{(m)}\in E\}.
$$
It is clear that $\hat{T}^{-1}E^{(m)} = E^{(m-1)}, m\geq 1 .$
\par If $\eta=\{E_{1},\ldots,E_{t}\}$ is a measurable partition of $\Omega$ then
we denote by $\eta^{(m)}$ the measurable partition of
$\hat{\Omega}$ defined by
$$\eta^{(m)} = \{E_{1}^{(m)},\ldots,E_{t}^{(m)}\}.
$$

Let $\xi$ be the zero-time partition of $\Omega$; $ \xi =
\{C_{0}(0),\ldots,C_{0}(r-1)\}$ where\\
$C_{0}(i) = \{x\in \Omega; x_{0} = i\},\ i\in \mathbb{Z}_{r}.$ For
$i, j\in\Bbb Z, i\leq j$ we put $\xi(i,j) =
\bigvee_{u=i}^{j}\sigma^{-u}\xi.$
\par Note that the corresponding entropies on
$(\hat{\Omega},\hat{\textbf{M}},\hat{\mu},\hat{T})$ are coincide
$(\Omega,\textbf{M},\mu,T)$ (cf. [2], [8])

\section{Main Results}

\par Let $\xi$ be a zero-time partition of $\Omega$, i.e.
$\xi=\{C_{0}(i), \, i\in \mathbb{Z}_{r} \}$ and $\{(m_{i},n_{i})
\}$ be a sequence of the lattice $\mathbb{Z}\times\mathbb{N}$.
Define a sequence of partitions of space $\Omega$ with respect to
$\mathbb{Z}^{2}$-action $\Phi$ by formula
$$
  {\xi}_{(m_{i}, n_{i})}={\Phi}^{(m_{i}, n_{i})}\xi, \qquad
i=\overline{1, \infty}.
$$
\par \textbf{Lemma 3.1.} Let $\xi_{(m_{i},n_{i})}\searrow \zeta$ and $\eta$ be an
arbitrary measurable partition with $ H_{\mu} (\xi
_{(m_{i},n_{i})} \ |\ \eta )< \infty $. Then
$$
H_{\mu} (\xi _{(m_{i},n_{i})} \ |\ \eta )\searrow H_{\mu}(\zeta \
|\ \eta ).
$$
\begin{proof}
Put $\alpha(n)=a+{\omega}^{-1}(n+1)-[a+{\omega}^{-1}]$, where
$[a]$ denotes the greatest integer $\leq a$. Denote
$$ \eta =\bigvee\limits_{{\alpha}(0)\le m_{i} \le {\alpha}(0)+r}
\Phi ^{(m_{i},1)}\xi
$$
Let $\xi_{(m_{i}, n_{i})}$ and $\zeta$ be two partitions as
$$
\xi _{(m_{i},n_{i})}=\bigvee\limits_{{\alpha}(0)\le m_{i} \le
{\alpha}(0)+r+2s} \Phi ^{(m_{i},0)}\xi\vee
\bigvee\limits_{n_{i}<0} \, \bigvee\limits_{{\alpha}(n_i)\le m_{i}
\le {\alpha}(n_i)+2s} \Phi ^{(m_{i},n_{i})}\xi
$$
and
$$
\zeta=\bigvee\limits_{n_{i}<0} \, \bigvee\limits_{{\alpha}(n_i)\le
m_{i}} \Phi ^{(m_{i},n_{i})}\xi.
$$

\par Denote $C_{\eta}(x)$, $C_{\zeta}(x)$ and $C_{\xi
_{(m_{i},n_{i}})}(x)$ elements of partitions $\eta$, $\zeta$ and
$\xi _{(m_{i},n_{i})}$ containing $x\in \Omega$, respectively.
Using Doob's theorem on convergence of conditional probabilities,
we have if $\xi_{(m_{i}, n_{i})}\searrow \zeta$ then
$$
\mu(C_{\xi_{(m_{i},n_{i})}}(x)|C_{\eta}(x))\rightarrow
\mu(C_{\zeta}(x)|C_{\eta}(x)).
$$
From this immediately follows that
$$
\lim\limits_{i \to \infty} \mu(\xi _{(m_{i},n_{i})} \ |\ \eta
)=\mu (\zeta \ |\ \eta ).
$$
From this using the properties of continuity of conditional
entropy  and logarithm we obtain that
$$
\lim\limits_{i \to \infty}H_{\mu}(\xi _{(m_{i},n_{i})} |\ \eta
)=H_{\mu}(\zeta \ |\ \eta ).
$$
\end{proof}
\par Now define a transformation $Q$ in the space of segments $I(a,
\omega)$ by
$$
  Q(I(a,\omega))=I(a', \omega),
$$
where $a'=a+{\omega}^{-1}$. Using properties of the
measure-theoretical entropy of dynamical system we shall prove the following theorem.\\

\par \textbf{Theorem 3.2.} If  $Q^{i} ( I_{1} )\subset Q^{i} (I)$, then
$H(Q^{i} ( I_{1} ))\leq H(Q^{i} (I))$.

\begin{proof} Let $(\Omega, M, \mu, \sigma)$ be
symbolic dynamic system and $\Phi^{(p, q)}=
\sigma^{p}T_{F_{[-k,k]}}^{q}$ be a $\mathbf{Z}^{2}$-action on
product space $\Omega$. Let $Q^{i}(I)$ and $Q^{i}(I_1)$, $i \ge
0$, be two transformations in the space of segments $I$ and
$I_{1}$, respectively. We consider the case when $i=0$. Other
cases can be shown in the same way. We have $H(I_{1})=
H_{r}(I_{1})+ H_{l}(I_{1})$. Since $\xi$ is a partition of
$\Omega$ we get
$$
\eta =\bigvee\limits_{a+{\omega}^{-1} \le p}
\sigma^{p}T_{F[-k,k]}^{b}\xi\preceq \bigvee\limits_{a+\omega
^{-1}\leq p} \sigma^{p}T^{1}_{F[-k,k]} \xi =
\bigvee\limits_{a+\omega ^{-1}\leq p} T_{F[-k+p,k+p]}^{1}\xi $$
From the continuity of conditional entropy and from Lemma 3.1 it
follows
$$
H (\bigvee\limits_{a+{\omega}^{-1} \le p}\sigma^{p}
T_{F[-k,k]}^{b}\xi \, {\vert} \, \bigvee\limits_{q=0}^{\infty} \,
\bigvee\limits_{|pm-(a+\omega ^{-1}q)|\leq s}
\sigma^{p}T_{F[-k,k]}^{-q}\xi) \le
$$
$$
\le H(\bigvee\limits_{a+{\omega}^{-1} \le p}\sigma^{p}
T_{F_{[-k,k]}}^{1}\xi \,  {\vert} \, \bigvee\limits_{q=0}^{\infty}
\, \bigvee\limits_{|p-(a+\omega^{-1}q)|\leq s}
\sigma^{p}T_{F[-k,k]}^{-q}\xi)
$$
It means that $H_{r} ( I_{1} ) \leq H_{r} ( I )$.
\par Similarly,
it can be shown that $H_{l} ( I_{1} ) \leq H_{l} ( I )$. From this
and the fact that $Q^{0}(I_{1})=I_{1}$, $Q^{0}(I)=I$ it is easily
follows the assertion of theorem 3.2 for the case $i=0$.
\end{proof}

\par Here, we investigate the measure-theoretic entropy of u-th
iteration of additive one-dimensional cellular automata. Recall
that the CA-map $T_{F[-k, k]}$ preserves the Bernoulli measure and
is non-invertible map of $\Omega$ generated
by a block map. So we should consider the condition $u\geq0$.\\

\par \textbf{Theorem 3.3.} Let $T_{F[-k,k]}$ be additive one-dimensional
cellular automata. Then for every $u\ge 0$ we have
$$
h_{\mu }T_{F[-k,k]}^{u})=uh_{\mu }(T_{F[-k,k]}).
$$

\begin{proof} Define the cylinder set $_{s}[i_{s},\ldots, i_{t}]_{t}=
\{x \in \Omega: \, x_{j}= i_{j}, \, s \le j \le t, \ i_{j}\in
\mathbb{Z}_{r} \}. $ Using the definition of partition $\xi (-k,
k)$ it can be easily checked
that\\
$\xi(-k,k)=\{_{-k}[i_{-k},...,i_{k}]_{k}:i_{j} \in
\mathbb{Z}_{r}\}$. Moreover, the partition $\xi (-k, k)$ is a
generator for $T_{F[-k,k]}$, i.e.

$$
\bigvee\limits_{i=0}^{\infty}T_{F[-k,k]}^{-i}\xi(-k,k)
=\varepsilon
$$
Using the properties of the measure-theoretic entropy and
Kolmogorov-Sinai theorem (cf. [9]) we get
\begin{eqnarray*} h_{\mu }(T_{F_{[-k,k]}}^{u})& = &h_{\mu}(T_{F[-k,k]}^{u},
\bigvee\limits_{i=0}^{u-1}T_{F[-k,k]}^{-i}\xi(-k,k)\\
& = &\lim\limits_{n \to \infty}\frac{1}{n}H_{\mu}
\left(\bigvee\limits_{j=0}^{n-1}T_{F[-k,k]}^{-uj}
(\bigvee\limits_{i=0}^{u-1}T_{F[-k,k]}^{-i}\xi(-k,k))\right)\\
& = & \lim\limits_{n \to \infty}\frac{u}{nu}H_{\mu} (
\bigvee\limits_{i=0}^{un-1}T_{F[-k,k]}^{-i}\xi(-k,k))\\
& = &u.2k\log r=uh_{\mu }(T_{F[-k,k]},\xi (-k,k)) = uh_{\mu
}(T_{F_{[-k,k]}}).
\end{eqnarray*}
\end{proof}
\par \textbf{Theorem 3.4.}
Let $\Phi =\sigma T_{F[-k,k]}$ be
$\mathbb{Z}\times\mathbb{N}$-action. Then for all $u\ge 0$,
$h_{\mu }(\Phi^{u})=uh_{\mu }(\Phi)$ and if the automaton rule
$F[-k, k]$ is bipermutative then
$h_{\vec{v}}(\Phi^{u})=uh_{\vec{v}}(\Phi)$ for all ${\vec{v}}\in
\mathbb{Z}^{2}$.
\begin{proof} Again first it is easy to see that the partition
$\xi(-k,k)=\{_{-k}[i_{-k},...,i_{k}]_{k}:i_{j} \in
\mathbb{Z}_{r}\}$ is generator for $\Phi =\sigma T_{F[-k,k]}$ that
is $ \bigvee\limits_{i=0}^{\infty} \sigma^{-i} T_{F[-k,k]}^{-i}\xi
(-k,k)=\varepsilon.$ So we have
\begin{eqnarray*}h_{\mu}(\Phi ^{u})& = &h_{\mu }(\Phi^{u},
\bigvee\limits_{i=0}^{u-1}{\Phi}^{-i} \xi(-k,k))\\
& = &\lim\limits_{s \to \infty}\frac{1}{s}H_{\mu}
\left(\bigvee\limits_{j=0}^{s-1}T_{F[-k,k]}^{-uj}{\sigma}^{-uj}
(\bigvee\limits_{i=0}^{u-1}T_{F[-k,k]}^{-i}{\sigma}^{-i}\xi(-k,k))\right)\\
& = & \lim\limits_{s \to \infty}\frac{1}{s}H_{\mu}
\left(\bigvee\limits_{j=0}^{s-1}T_{F[-k,k]}^{-uj}{\sigma}^{-uj}
(\bigvee\limits_{i=0}^{u-1}T_{F[-k,k]}^{-i}\xi(-k-i,k-i))\right)\\
& = &\lim\limits_{s \to \infty}\frac{1}{s}H_{\mu}
\left(\bigvee\limits_{j=0}^{s-1}T_{F[-k,k]}^{-uj}
(\bigvee\limits_{i=0}^{u-1}T_{F[-k,k]}^{-i}\xi(-k-(i+ju),
k-(i+ju)))\right)\\
 & = &u \lim\limits_{s \to
\infty}\frac{1}{us}H_{\mu}
\left(\bigvee\limits_{j=0}^{us-1}T_{F[-k,k]}^{-j} \xi(-k-(i+ju),
k-(i+ju))\right)=uh_{\mu}(\Phi)
\end{eqnarray*}

\par Now we consider the directional entropy of
$\mathbb{Z}^{2}$-action. Here we only consider ${\vec{v}}\in
\mathbb{Z}^{2}$. Using the definition of $h_{\vec{v}}(\Phi)$ (cf.
[2]) we have
\begin{eqnarray*} h_{\vec{v}}(\Phi^{u}) & =
&h_{\hat{\mu }}(\hat{\sigma
}^{up}\hat{T}_{F[-k,k]}^{uq})\\
& = & h_{\mu }(\sigma^{up}T_{F[-k,k]}^{uq})\\
& = & h_{\mu }(\sigma^{up}T_{F^{uq}[-qku,qku]})\\
 & = &h_{\mu}(T_{F^{uq}[-qku+up,qku+up]})=uh_{\vec{v}}(\Phi).
\end{eqnarray*}

\end{proof}
\par One can investigate for ${\vec{v}}\in \mathbb{R}^{2}$.

\vspace{1cm}

\baselineskip=0.4\baselineskip

{\small
\noindent Harran University \\
Arts and Sciences Faculty \\
Department of Mathematics \\
 63100, \c Sanl\i urfa, TURKEY \\
e-mail address : akinhasan@harran.edu.tr

\end{document}